\documentclass[11pt]{amsart}
\usepackage{amsmath}
\usepackage{amsfonts}
\usepackage{amssymb}
\newtheorem{thm}{Theorem}

\newtheorem{cor}[thm]{Corollary}

\newtheorem{lemma}[thm]{Lemma}
\newtheorem{prop}[thm]{Proposition}

\newtheorem{defn}[thm]{Definition}

\def \H{{\mathbb H}}
\newcommand{\bb}[1]{\mathbb{#1}}
\newcommand{\cl}[1]{\mathcal{#1}}

\begin{document}

\title[Kadison-Singer]{Projections and the Kadison-Singer Problem}

\author[P. Casazza]{Pete Casazza}
\address{Department of Mathematics, University of Missouri}
\email{pete@math.missouri.edu}

\author[D. Edidin]{Dan Edidin}
\address{Department of Mathematics, University of Missouri}
\email{edidin@math.missouri.edu}

\author[D. Kalra]{Deepti Kalra}
\address{Department of Mathematics, University of Houston, Houston, Texas 77204-3476, U.S.A.}
\email{deepti@math.uh.edu}

\author[V.~I.~Paulsen]{Vern I.~Paulsen}
\address{Department of Mathematics, University of Houston,
Houston, Texas 77204-3476, U.S.A.}
\email{vern@math.uh.edu}

\thanks{The first author was supported by NSF DMS 0405376 and the
  third and fourth authors were supported by NSF DMS 0600191.}
\subjclass[2000]{Primary 46L15; Secondary 47L25}

\begin{abstract}
We prove some new equivalences of the paving conjecture and obtain
some estimates on the paving constants. In addition we give a new
family of counterexamples to one of the Akemann-Anderson conjectures. 
\end{abstract}

\maketitle

\section{Introduction}

Let $\cl H$ be a separable, infinite dimensional Hilbert 
space and let $B(\cl H)$ denote the bounded, linear operators 
on $\cl H.$ By a  {\em MASA} we mean a maximal, abelian subalgebra 
of $B(\cl H)$. 
R. Kadison and I. Singer studied \cite{KS} whether or not pure 
states on a MASA extend uniquely to states on $B(\cl H).$ In their 
original work on this subject \cite{KS}, it was shown that this 
question has a negative answer if the MASA had any continuous part. The remaining case, whether or not
pure states on discrete MASA's have unique extensions to states 
on $B(\cl H)$ has come to be known as the {\em Kadison-Singer problem.} 
their work showed that this problem was equivalent to certain questions 
about "paving" operators by projections. J. Anderson\cite{A} developed 
this idea significantly into a series of so-called "paving" conjectures. 
Since that time there has been a great deal of research on these paving 
conjectures \cite{AA}, \cite{BHKW}, \cite{BHKW2}, \cite{BT}, \cite{C}, 
\cite{CT}, \cite{HKW}, \cite{HKW2} and \cite{HKW3}.

In this paper, we begin by restating some of these paving conjectures and add a few new equivalent paving conjectures.


\section{Some New Equivalences of the Paving Conjecture}

Let us begin with the familiar.

Given $A \subseteq I,$ where $I$ is some index set, we let 
$Q_A \in B(\ell^2(I))$ denote the
diagonal projection defined by $Q_A= ( q_{i,j}), q_{i,i}=1, i \in A,
q_{i,i}= 0, i \notin A$ and $q_{i,j} =0, i \ne j.$

\begin{defn} An operator $T \in B(\ell^2(I))$ is said to have an
  $(r,\epsilon)$-paving if there is a partition of $I$ into $r$
  subsets $\{ A_j \}^r_{j=1}$ such that $\|Q_{A_{j}} T Q_{A_j} \| \le
  \epsilon.$ A collection of operators $\cl C$ is said to be $(r,
  \epsilon)$-pavable if each element of $\cl C$ has an $(r,
  \epsilon)$-paving.
\end{defn}

Note that in this definition, we do {\em not} require that the
diagonal entries of the operator be 0.

Some classes that will play a role are:
\begin{itemize}
\item $\cl C_{\infty} = \{ T=(t_{i,j}) \in B(\ell^2(\bb N)): \|T\| \le 1,
t_{i,i}=0 \forall i \in \bb N \},$
\item $\cl C= \cup_{n=2}^{\infty} \{ T=(t_{i,j}) \in M_n: \|T\| \le 1,
  t_{i,i} = 0, i=1,...,n \},$
\item $\cl S_{\infty} = \{ T \in \cl C_{\infty}: T=T^* \},$
\item $\cl S= \{ T \in \cl C : T=T^* \},$
\item $\cl R_{\infty}= \{T \in \cl S_{\infty}: T^2=I \},$
\item $\cl R = \{T \in \cl S : T^2=I \},$
\item $\cl P^{\infty}_{1/2} = \{ T=(t_{i,j}) \in B(\ell^2(\bb N)):
  T=T^*=T^2, t_{i,i} =1/2, \forall i \in \bb N \},$
\item $\cl P_{1/2} = \cup_{n=2}^{\infty} \{ T=(t_{i,j}) \in M_n :
  T=T^*=T^2, t_{i,i} =1/2, i=1,...,n \}.$
\end{itemize}

Note that the operators satisfying, $R=R^*, R^2=I$ are {\em reflections} 
and that for such an operator, $\sigma(R) = \{ -1, +1 \}$. Since, the 
traces of our matrices are 0, in the finite dimensional case, these 
types of reflections can only exist in even dimensions. If the space 
is $2n$-dimensional, then there exists an $n$-dimensional subspace 
that is fixed by $R$ and such that for any vector $x$ orthogonal to 
the subspace. $Rx= -x$. 

J. Anderson's \cite{A} remarkable contribution follows.

\begin{thm}[Anderson] The following are equivalent:
\begin{enumerate}
\item the Kadison-Singer conjecture is true,
\item for each $T \in \cl C_{\infty},$ there exists
  $(r,\epsilon)$(depending on T) $\epsilon < 1$, such that 
$T$ is $(r, \epsilon)$-pavable,
\item there exists $(r, \epsilon), \epsilon <1,$ such that 
$\cl C_{\infty}$ is
  $(r,\epsilon)$-pavablle,
\item there exists $(r,\epsilon), \epsilon < 1$, such that $\cl C$ is
  $(r,\epsilon)$-pavable,
\item for each $T \in \cl S_{\infty},$ there exists $(r,\epsilon),
  \epsilon < 1$(depending on T), such that $T$ is $(r,
  \epsilon)$-pavable,
\item there exists $(r, \epsilon), \epsilon < 1,$ such that $\cl
  S_{\infty}$ is $(r, \epsilon)$-pavable,
\item there exists $(r, \epsilon), \epsilon < 1,$ such that $\cl S$ is
  $(r, \epsilon)$-pavable.
\end{enumerate}
\end{thm}

Generally, when people talk about the paving conjecture they mean one
of the above equivalences of the Kadison-Singer problem. Also,
generally, when one looks at operators on an infinite dimensional
space, it is enough to find $(r,\epsilon)$ depending on the operator,
but for operators on finite dimensional spaces it is essential to have
a uniform $(r, \epsilon),$ for all operators of norm one. Finally, since 
$\cl S_{\infty} \subset \cl C_{\infty}$,
people looking for counterexamples tend to study $\cl C_{\infty}$, while people
trying to prove the theorem is true, study $\cl S_{\infty}$ or $\cl S.$
However, by the above equivalences, if a counterexample exists in one
set then it must exist in the other as well.

In this spirit, we prove that the following smaller sets with ``more
structure'' are sufficient for paving.

\begin{thm}\label{Main} 
 Let $\epsilon < 1,$ then the following are equivalent:
\begin{enumerate}
\item the set $\cl S_{\infty}$ can be $(r_1, \epsilon)$-paved,
\item the set $\cl R_{\infty}$ can be $(r_1, \epsilon)$-paved,
\item the set $\cl P^{\infty}_{1/2}$ can be 
$(r_2, \frac{1+ \epsilon}{2})$-paved,
\item the set $\cl S$ can be $(r_1, \epsilon)$-paved,
\item the set $\cl R$ can be $(r_1, \epsilon)$-paved,
\item the set $\cl P_{1/2}$ can be $(r_2, \frac{1+ \epsilon}{2})$-paved.
\end{enumerate}
\end{thm}

\begin{proof} Since the reflections are a subset of the self-adjoint 
matrices, it is clear that (1) implies (2) and that (4) implies (5).

To see that (2) implies (1), let $A \in \cl S_{\infty},$ and set 
$$R= \begin{pmatrix} A & \sqrt{I-A^2}\\ \sqrt{I- A^2} & -A \end{pmatrix},$$
then $R \in \cl R_{\infty}$ and clearly any $(r, \epsilon)$-paving of $R$ 
yields an $(r, \epsilon)$-paving of $A$.

Thus, (1) and (2) are equivalent and similarly, (4) and (5) are equivalent.

To see the equivalence of (2) and (3), note that $R \in \cl R_{\infty}$
 (respectively, $\cl R$) if and only if 
$P= (I+R)/2 \in \cl P^{\infty}_{1/2}$ (respectively, 
$\cl P_{1/2}$). Also, if $\|Q_AR Q_A \| \le \epsilon,$ 
then $\|Q_A P Q_A \| \le (1+ \epsilon)/2.$ Thus, if 
$\cl R_{\infty}$ can be $(r_1, \epsilon)$-paved, 
then $\cl P^{\infty}_{1/2}$ can be $(r_1, \frac{1 + \epsilon}{2})$-paved.

Conversely, given $R \in \cl R_{\infty}$, let $P= (I+R)/2.$ If 
$\| Q_A P Q_A \| \le (1+ \epsilon)/2= \beta,$ then, 
$$0 \le Q_A PQ_A \le \beta Q_A,$$ and since $R= 2P - I,$ 
we have that $$-Q_A \le Q_ARQ_A \le (2 \beta -1) Q_A = \epsilon Q_A.$$ 
Applying the same reasoning to the reflection $-R$, we get a new 
projection, $P_1 = (I - R)/2$, with a possibly different paving 
of $P_1$, such that $-Q_B \le Q_B(-R)Q_B \le \epsilon Q_B.$ 
Thus, $-\epsilon Q_B \le Q_BRQ_B$ and if $Q_C= Q_AQ_B,$ 
we have that $-\epsilon Q_C \le Q_CRQ_C \le +\epsilon Q_C.$  
Therefore, we have that the set of all products of the $Q_A$'s 
and $Q_B$'s pave $R.$  Thus, if $\cl P^{\infty}_{1/2}$ can 
be $(r_2, \frac{1+ \epsilon}{2})$-paved, then $\cl R_{\infty}$ 
can be $(r_2^2, \epsilon)$-paved.

The proof of the equivalence of (5) and (6), is identical.

Finally, (1) and (4) are equivalent by the standard limiting argument. In particular, see \cite[Proposition~2.2]{C} and the proof of \cite[Theorem~2.3]{C}.

\end{proof}
 
\begin{cor}  The following are equivalent:
\begin{enumerate}
\item the Kadison-Singer conjecture is true,
\item for each $R \in \cl R_{\infty}$ there is a $(r, \epsilon),
  \epsilon < 1$(depending on R) such that $R$ can be $(r,
  \epsilon)$-paved,
\item there exists $(r,\epsilon), \epsilon < 1,$ such that every $R
  \in \cl R$ can be $(r,\epsilon)$-paved,
\item for each $P \in \cl P^{\infty}_{1/2}$ there is a $(r,\epsilon),
  \epsilon < 1$(depending on P) such that $P$ can be
  $(r,\epsilon)$-paved,
\item there exists $(r,\epsilon), \epsilon < 1,$ such that every $P
  \in \cl P_{1/2}$ can be $(r, \epsilon)$-paved.
\end{enumerate}
\end{cor}

We will need some results from frame theory 
in this paper.  We refer the reader to \cite{Ch} for these.
We will briefly give the definitions we will be using.  If 
$\{f_i\}_{i\in I}$ is a family of vectors in a Hilbert space
$\H$, the {\em analysis operator} of this family is $T:\H
\rightarrow \ell_2(I)$ given by $T(f) = \{\langle f,f_i\rangle\}_{i\in I}$,
and the {\em synthesis operator} is $T^{*}(\{a_i\}_{i\in I}) =
\sum_{i\in I}a_if_i$.  If $T$ is bounded, we call $\{f_i\}_{i\in I}$
a {\em Bessel sequence}, if it is also onto we call
this a {\em frame}, and if $T$ is invertable it is
a {\em Riesz basis}.  The frame is {\em equal-norm} or
{\em uniform} if the $f_i$ all have
the same norm and it is {\em equiangular} if there is a constant
$c$ so that $|\langle f_i,f_j\rangle|=c$ for all $i\not= j\in I$.
This is a {\em Parseval frame} if $T$ is a partial isometry.  In this
case, the Gram matrix $(\langle f_i,f_j\rangle)_{i,j\in I}$ is
an orthogonal projection of $\ell_2(I)$ onto the range of the
analysis operator (and this takes $e_i$ to
$T(f_i)$ where $\{e_i\}_{i\in I}$ is the unit vector basis
of $\ell_2(I)$).  If dim $\H = k$ and $|I|=n$ we call this an
{\em (n,k)-frame}.

A sort of meta-corollary or Theorem \ref{Main} is that the 
frame based conjectures that are known to be equivalent to the
Kadison-Singer result can be reduced to the case of uniform Parseval
frames of redundancy 2. Similarly, for most harmonic analysis
analogues of paving, it is enough to consider say subsets $E \subseteq
[0,1]$ of Lebesgue measure 1/2.
We state one such equivalence.  The {\em Feichtinger Conjecture} in
frame theory asserts that every unit norm Bessel sequence is a finite
union of Riesz basic sequences.  
Casazza and Tremain \cite{CT} have shown that
the Feichtinger conjecture is equivalent to Kadison-Singer.

\begin{thm} The Feichtinger conjecture is true if and only if for each
  Parseval frame $\{ f_n \}_{n \in \bb N}$ for a Hilbert space with
  $\|f_n\|^2=1/2$ $ \forall n$ there is a partition $\{ A_j \}_{j=1}^r$
  of $\bb N$ into r disjoint subsets (with r depending on the frame)
  such that for each $j, \{ f_n \}_{n \in A_j}$ is a Riesz basis for
  the space that it spans.
\end{thm}  
\begin{proof} Clearly, if the Feichtinger is true, then it is true for
  this special class of frames.

Conversely, assume that the above holds and let $P \in \cl
P^{\infty}_{1/2}$. Then there exists a Parseval frame $\{ f_n \}_{n
  \in \bb N}$ for some Hilbert space $\cl H,$ such
that $I-P= ( \langle f_j, f_i \rangle)$ is their Grammian.
Now let $\{ A_k \}_{k=1}^r$ be the partition of $\bb N$ into r
disjoint subsets as above and let 
$\cl H_k = \overline{span} \{ f_n : n \in A_k \}$ denote 
the closed linear span.

Since $\{ f_n : n \in A_k \}$ is a Riesz basis for $\cl H_k$, there
exists an orthonormal basis, $\{ e_n: n \in A_k \}$ for $\cl H_k$ and
a bounded invertible operator, $S_k: \cl H_k \to \cl H_k,$ with
$S_k(e_n) = f_n.$

We have that $Q_{A_k}(I-P)Q_{A_k} = ( \langle f_j,f_j \rangle)_{i,j \in A_k} = (\langle S_k^*S_ke_j,
e_i \rangle ) \ge c_k Q_{A_k}$ where $S_k^* S_k \ge c_k Q_{A_k}$ for
some constant $0 <c_k \le 1$ since $S_k$ is invertible. Hence,
$Q_{A_k}PQ_{A_k} \le (1-c_k)Q_{A_k}$ and we have that, $max \{
\|Q_{A_k}PQ_{A_k} \| : 1 \le k \le r \} < 1.$

Hence, condition (5) of Corollary~4 is met and so Kadison-Singer is
true and thus, by \cite[Theorem~5.3]{CT}, the Feichtinger conjecture is
true.
\end{proof}


\section{Some Paving Estimates}

In this section we derive some estimates on paving constants that give
some basic relationships between $r$ and $\epsilon.$ In
particular, we will prove that $\cl P_{1/2}$ cannot be $(2,
\epsilon)$-paved for any $\epsilon <1.$

We begin with a result on paving $\cl R.$

\begin{thm} Assume that $\cl R$ is
  $(r, \epsilon)$-pavable. Then $1 \le r \epsilon^2.$
\end{thm}
\begin{proof}  Recall that an $n \times n$ matrix C is a {\em conference matrix} if $C=C^*,
  c_{i,i}=0, c_{i,j}= \pm 1, i \ne j$ and $C^2= (n-1)I.$ Such matrices
  exist for infinitely many $n$.

Set $A= \frac{1}{\sqrt{n-1}}C,$ then $A$ is a unitary matrix with
zero diagonal.

Assume that $\{1,...,n \} = B_1 \cup ... \cup B_r$ is a partition such
that $\|Q_{B_i}AQ_{B_i} \| \le \epsilon.$
Let $d = max \{ card(B_i) \}$ and let $B_j$ attain this max. Note that
$d \ge \frac{n}{r}.$
Set $A_j = Q_{B_j}AQ_{B_j},$ then the Schur product $A_j*A_j=
\frac{1}{n-1}[J_d- I_d]$ where $J_d$ denotes the matrix of all 1's.
Hence, $\frac{d-1}{n-1} = \|A_j*A_j \| \le \|A_j\|^2 \le \epsilon^2.$
Thus, $\frac{n/r - 1}{n-1} \le \epsilon^2,$ and the result follows 
by letting $n \to +\infty$
\end{proof}

\begin{prop}\label{AP201}
If every projection $P \in \cl P_{1/2}$ can be
$(r, \epsilon)$-paved then every projection $Q$ with
$$
\frac{1}{2}-\delta \le \langle Qe_i,e_i\rangle \le \frac{1}{2}
+\delta,
$$
can be $(r,\beta)$-paved, where 
$$
\beta = (1+2\delta)\epsilon,
$$
and so $\beta < 1$ when $\delta$ is small enough.
\end{prop}

\begin{proof}
Let $D$ be the diagonal of $Q$ and let $B=Q-D$.  Then
$$
\|B\|\le \frac{1+2\delta}{2}.
$$
To see this note that for any vector $x$
$$
0\le \langle Bx,x\rangle + \langle Dx,x\rangle \le 1,
$$
since $Q$ is a projection.  Hence,
$$
-\langle Dx,x\rangle \le \langle Bx,x\rangle \le
1-\langle Dx,x\rangle.
$$
Hence,
$$
\|B\| = \sup_{\|x\|=1}|\langle Bx,x\rangle | \le max\{|\langle Dx,x\rangle |,
|1-\langle Dx,x\rangle |\} \le \frac{1}{2} + \delta= \frac{1+2 \delta}{2}.
$$

Let $R=R^{*}$ be the symmetry we get by dilating
$$
\frac{2}{1+2\delta}B,
$$
as in the proof of Theorem~3.  Let $P = \frac{1}{2}(I+R)$ be the
projection with $1/2's$ on the diagonal.  If we can $(r,\epsilon)$-pave
$P$ with $\{A_j\}_{j=1}^{r}$ then we have $(r,\epsilon)$-paved
$$
\frac{1}{2}I+\frac{1}{1+2\delta}B.
$$
Substituting $B=Q-D$ we have an $(r,\epsilon)$-paving of
$$
\frac{1}{1+2\delta}Q+\frac{1}{2}I-\frac{1}{1+2\delta}D =
\frac{1}{1+2\delta}\left ( Q+\frac{1+2\delta}{2}I-D \right ).
$$
Now, for any $j=1,2,\ldots ,r$ since 
$$
\frac{1+2\delta}{2}I-D
$$ 
is a positive operator,
$$
\|Q_{A_j}QQ_{A_j}\| \le \|Q_{A_j} (Q+\frac{1+2\delta}{2}I-D)Q_{A_j}\|
\le (1+2\delta)\epsilon < 1.
$$
\end{proof}

\begin{thm}
Assume that $\cl P_{1/2}$ can be 
$(r, \epsilon)$-paved. Then $\frac{r}{2(r-1)}\le \epsilon.$

Note: When r=2 this implies that $\epsilon=1,$ and hence 2-paving is 
impossible.
\end{thm}

\begin{proof}
Let $m > 2$ be an integer and consider a uniform, Parseval 
(n,k)-frame with $n= mr, k= m(r-1)+1$. This will give rise to a 
projection $Q$ with 
diagonal 
entries, $\frac{m(r-1)+1}{mr}= \frac{1}{2} + \delta,$ where 
$\delta=\frac{m(r-2)+3}{2mr}.$  To see this, let
\begin{eqnarray*}
\delta &=& \frac{m(r-1)+1}{mr}- \frac{1}{2}\\
&=& \frac{2[m(r-1)+1] - mr}{2mr}\\
&=& \frac{m(r-1)+2}{2mr}.
\end{eqnarray*} 

By the above result, $Q$ can be $(r, \beta)$-paved, where $\beta=(1+2 
\delta) \epsilon.$

However, for any r paving of Q, one of the blocks must be of size at least 

$$
n/r = m = n-k+1,
$$
 by the choice of $n$ and $k$. Since $Q$ is a rank $k$ projection, 
this block will have 
norm 1 by the eigenvalue inclusion principle or by the eigenvalue
 interlacing results. Hence $\beta \ge 1.$  We solve for $\epsilon$:
$$
(1+2\delta) \epsilon \ge 1.
$$
So
\begin{eqnarray*}
\epsilon &\ge& \frac{1}{1+2\delta}\\
&=& \frac{mr}{m(2r-2) +2}.
\end{eqnarray*}
Letting $m \to + \infty$ yields 
$$
\epsilon \ge \frac{r}{2(r-1)}.
$$ 
\end{proof}

\begin{cor}
The set $\cl P_{1/2}$ is not
$2$-pavable.
\end{cor}

\begin{cor}
The set $\cl R$ is not $2$-pavable.
\end{cor}

We now generalize the results of the last theorem.

\begin{thm}
For each $r,n \in \bb N$ with $r>1$ there is an $\epsilon_n>0$ so 
that whenever
$P$ is a projection on $\ell_2^n$ with 
$\frac{1}{r} \le \langle Pe_i,e_i \rangle \le 1-\frac{1}{r}$ for
all $i=1,2,\ldots ,n$ then $P$ is $(r,1-\epsilon_n)$-pavable.

Moreover, for any $\delta >0$ there is an $n\in \bb N$ and
a projection $P$ on $\ell_2^{2n}$ of rank $n$ so that
$\frac{1}{r} - \delta \le \langle Pe_i,e_i \rangle \le 1-\frac{1}{r}
+ \delta$ for all $i=1,2,\ldots ,2n$ while $P$ is not 
$(r,\epsilon)$-pavable for any $\epsilon <1$.
\end{thm}

\begin{proof}
Given our assumptions, we will check the Rado-Horn Theorem 
(see \cite{CKS} and its references) to see
that the row vectors of our projection can be divided into
$r$ linearly independent sets.  Then the rest of the first part
of the theorem follows by the same argument (adjusted for $r$)
as in the last theorem.  For any $J\subset \{1,2,\ldots ,2n\}$ let
$P_J$ be the orthogonal projection of $\ell_2^{2n}$ onto the
span $\{Pe_i\}_{i\in J}$.  Now,
$$
dim\ span\ \{Pe_i\}_{i\in J} = \sum_{i=1}^{2n}\|P_J Pe_i\|^2
\ge \sum_{i\in J}\|Pe_i\|^2 \ge |J| \frac{1}{r}.
$$
By the Rado-Horn Theorem we can now write $\{Pe_i\}_{i=1}^{2n}$
as a union of $r$-linearly independent sets.

For the moreover part, choose a $k\in \bb N$ so that
$$
\frac{1}{r} - \delta < \frac{k}{rk+1} \le \frac{1}{r}\le 1-\frac{1}{r}.
$$
Now, choose an $n$ so that
$$
\frac{1}{r} - \delta \le \frac{n-rk}{2n-(rk+1)}\le 1-\frac{1}{r}
+ \delta.
$$
With $\{e_i\}_{i=1}^{n}$ the unit vectors in $\ell_2^n$ we can
choose an equal norm Parseval frame $\{f_i\}_{i=1}^{rk+1}$ for
$\{e_i\}_{i=1}^{k}$.  Next, choose an equal norm Parseval frame
$\{f_i\}_{i=rk+2}{2n}$ for $\{e_i\}_{i=k+1}^{n}$.  Now,
$$
\frac{1}{r}-\delta \le \|f_i\|^2 = \frac{k}{rk+1}\le \frac{1}{r}
\le 1-\frac{1}{r},
$$
and
$$
\frac{1}{r} - \delta \le \frac{n - rk}{2n-(rk+1)} \le 1-\frac{1}{r}
+ \delta.
$$
Taking the embedding of this Parseval frame with $2n$-elements for
$\ell_2^n$ into $\ell_2^{2n}$
we get a projection $P$ on $\ell_2^{2n}$ which has
rank $n$ and looks like
$$
\begin{bmatrix} \|f_1\|^2 & b_{1,2} & \ldots & b_{1,(rk+1)} & 0 & 0 
& \ldots & 0\\
b_{21} & \|f_1\|^2 & \ldots & b_{2,(rk+1)} & 0 & 0 & \ldots & 0\\
\vdots & \vdots & \ldots & \vdots & \vdots & \vdots & \ldots & \vdots\\
b_{(rk+1),1} & b_{(rk+1),2} & \ldots &\|f_1\|^2 & 0 & 0 & \ldots & 0\\
0 & 0 & \ldots & 0 & \|f_{rk+2}\|^2 & a_{(rk+2),(rk+3)} & \ldots & 
a_{(rk+2),2n}\\
0 & 0 & \ldots & 0 &  a_{(rk+3),(rk+2)} & \|f_{rk+2}\|^2 &  \ldots & 
a_{(rk+3),2n}\\
\vdots & \vdots & \ldots & \vdots & \vdots & \vdots & \ldots & \vdots\\
0 & 0 & \ldots & 0 & a_{(2n),(rk+2)} & a_{(2n),(rk+3)} & \ldots 
& \|f_{rk+2}\|^2 
\end{bmatrix}
$$
For this projection, for any $J\subset \{1,2,\ldots n\}$ with
$|J|>r$ the family $\{e_i\}_{i\in J}$ is linealy dependent and
so $PQ_A$ has a zero eigenvalue.  Hence, $(I-P)Q_A$ has one as
an eigenvalue and hence is not $\epsilon$-pavable for any
$\epsilon >0$. 
\end{proof}


\section{Counterexamples to the Akemann-Anderson Conjecture}

In \cite{AA} Akemann and Anderson introduce two paving
conjectures, denoted {\em Conjecture A} and {\em Conjecture B.} 
They prove that Conjecture A implies Conjecture B and that Conjecture B
implies Kadison-Singer, but it is not known if either of these
implications can be reversed. Weaver\cite{W} provides a set of
counterexamples to Conjecture A.  Thus, if these three statements were
all equivalent then Weaver's counterexample would be the end of the
story. However, it is generally believed that Conjecture A is strictly
stronger than the Kadison-Singer conjecture.

In this section, we show that the Grammian projection matrices of {\em any}
uniform, equiangular (n,k)-frame, with $n>5k$ yield
counterexamples to Conjecture A. It is known that infinitely many such frames 
exist for arbitrarily large n and k. The significance of our new set of 
counterexamples is that by the results of J. Bourgain and L. 
Tzafriri \cite{BT}, there exists $\epsilon <1,$ such that the 
family of self-adjoint, norm one, 0 diagonal matrices obtained from 
these frames is $(2, \epsilon)$-pavable.  

Thus, these new examples drive an additional wedge between Conjecture A and 
Kadison-Singer.

We then turn our techniques to Conjecture B and derive some results
that could lead to a counterexample to Conjecture B.

We now describe the Akemann-Anderson conjectures.  
Let $P=(p_{i,j}) \in M_n$ be the matrix of a projection 
and set $\delta_P= max \{ p_{i,i}: 1 \le i \le n \}.$ By a 
{\em diagonal symmetry} we mean a diagonal matrix whose diagonal 
entries are $\pm 1,$ that is, S is a diagonal self-adjoint unitary.

{\bf Conjecture A \cite[7.1.1]{AA}.} For any projection, P there 
exists a diagonal symmetry, S, such that $\|PSP \| \le 2 \delta_P.$

{\bf Conjecture B \cite[7.1.3]{AA}.} There exists $\gamma, \epsilon > 0$ 
(and independent of n) such that for any P with $\delta_P < \gamma$ 
there exists a diagonal symmetry, S, such that $\|PSP\| < 1- \epsilon .$

Weaver\cite{W} states that a counterexample to Conjecture B would
probably lead to a negative solution to Kadison-Singer. We believe
that these two conjectures are really more closely related to
2-pavings and this is why we believe that counterexamples to
Conjecture B should be close at hand.

Finally, note that Conjecture B is about paving projections with small
diagonal. But our results show that Kadison-Singer is equivalent to
paving projections with diagonal 1/2. This would also seem to put
further distance between
these Akemann-Anderson conjectures
and the Kadison-Singer conjecture.

\begin{prop} Let $P= \begin{pmatrix} A & B\\B^* & C \end{pmatrix}$ be a 
projection, written in block-form with $A$ $m \times m, B$ $m \times (m+l), C$ $(m+l) \times (m+l)$, where $l \ge 0.$ Then there exists a $m \times m$ unitary $U_1$ and an $(m+l) \times (m+l)$ unitary $U_2$ such that, $U_1^*AU_1=D_1, U_1^*BU_2= (D_2,0), U_2^*CU_2= \begin{pmatrix} D_3 & 0\\0 & D_4 \end{pmatrix}$  where each of the $D_i$'s is a diagonal matrix with non-negative entries, $D_1, D_2, D_3$ are all $m \times m, D_4$ is $l \times l$ with $1$'s and $0$'s for its diagonal entries and the $0$'s represent matrices of all zeroes that are either $m \times l$ or $l \times m.$
\end{prop}
\begin{proof} First note that since $P$ is a projection we have that 
$A^2 + BB^*=A, B^*B +C^2 =C$ and $AB+BC=B.$ Also, since the rank of 
$B$ is at most $m$, the matrix $B^*B$ must have a kernel of dimension at 
least $l$.

Conjugating $P$ by a unitary of the form 
$U= \begin{pmatrix} I_m &0\\0 &U_2 \end{pmatrix}$, we may diagonalize 
$C$ and the new matrix, $P_1$, will still be a projection. Since 
$U_2^*B^*BU_2 = U_2^*(C - C^2)U_2,$ we see that both sides of this 
equation are in diagonal form. Since at least $l$ of the diagonal 
entries of $U_2^*B^*BU_2$ are zeroes, after applying a permutation 
if necessary, we may assume that, 
$$U_2^*B^*BU_2 = \begin{pmatrix} D_2^2 & 0\\0 & 0 \end{pmatrix}, U_2^*CU_2= 
\begin{pmatrix} D_3 &0\\0 & D_4 \end{pmatrix},$$
where $D_2, D_3, D_4$ are as claimed.

Now we may polar decompose the $m \times (m+l)$ 
matrix $BU_2= W |BU_2| = W \begin{pmatrix} D_2 & 0\\0 & 0 \end{pmatrix},$ 
where $W$ is a $m \times (m+l)$ partial isometry whose initial space is the 
range of $|BU_2|$. Thus, $W=(W_1, 0)$ where $W_1$ is an $m \times m$ 
partial isometry. Hence, we may extend $W_1$ to an $m \times m$ unitary 
$U_1$ with $W_1 D_2= U_1 D_2$ and 
$BU_2= (U_1,0) \begin{pmatrix} D_2 & 0\\0 & 0 \end{pmatrix}= (U_1D_2, 0).$

Conjugating $P_1$ by the unitary 
$\begin{pmatrix} U_1 & 0\\0 & I_{m+l} \end{pmatrix}$ we arrive at a 
new projection of the form,
$$ \begin{pmatrix} U_1^*AU_1  & D_2 & 0\\D_2 & D_3 & 0 \\0 & 0 & D_4 
\end{pmatrix}.$$

Note that since this last matrix is a projection, 
$U_1^*AU_1D_2 + D_2D_3 = D_2$ and so, $U_1^*AU_1D_2$ is diagonal. 
If all of the entries of $D_2$ were non-zero, then this would 
imply that $U_2^*AU_2$ is diagonal. In general, this implies that 
$U_1^*AU_1$ (which is self-adjoint) is of the form 
a diagonal matrix direct 
sum with another matrix corresponding to the block where $D_2$ is 0. 
Conjugating $U_1^*AU_1$ by another unitary to diagonalize this lower 
block, yields the desired form.

Finally, note that since $D_4$ is a diagonal projection, all of its 
entries must be $0$'s or $1$'s.
\end{proof}

\begin{lemma} Let $P= \begin{pmatrix} a & b\\b & c \end{pmatrix}$ be a 
non-zero projection  with real entries and let 
$S= \begin{pmatrix} 1 & 0\\0 & -1 \end{pmatrix}$. Then 
$\|PSP\| = |1 -2c|.$
\end{lemma}
\begin{proof} If $P$ is rank 2 then $P=I$ and the result is trivial. So 
assume that $P$ is rank one. We have that 
$PSP= \begin{pmatrix} a^2-b^2 & ab-bc\\ab-bc & b^2-c^2 \end{pmatrix}$ 
and since $P$ is a rank one projection, $a+c=1, b^2+c^2=c.$ A little 
calculation shows that the characteristic polynomial of $PSP$ is 
$x^2- Tr(PSP)x + Det(PSP) = x^2 - (1-2c)x,$ and hence 
the eigenvalues are 0 and 1-2c, from which the result follows.
\end{proof}

Note that when $S$ is a diagonal symmetry, then $\|PSP\|= \|P(-S)P\|$
and so we may and do assume in what follows that the number of $-1$'s
is greater than or equal to the number of $+1$'s. Also, given a matrix
$A$, we let $\sigma(A)$ denote the spectrum of $A$ and set $\sigma
^\prime(A) \equiv \sigma(A)\backslash \{0\}.$
 
\begin{thm} Let $P= \begin{pmatrix} A & B\\B^* & C \end{pmatrix}$
  be an $n \times n$ projection and let $S=\begin{pmatrix} I & 0\\0 &
    -I \end{pmatrix}$ be a diagonal symmetry. Then $\|PSP\| \ge \max
  \{| 1- 2\lambda|: \lambda \in \sigma^\prime(C) \cup \sigma^\prime(A) \}.$
\end{thm}
\begin{proof} Given any unitary of the type in the above Proposition,
  we have that $\|PSP\|= \|U^*PSPU\|=\|(U^*PU)(U^*SU)(U^*PU)\|=
  \|(U^*PU)S(U^*PU)\|.$ Thus, we may and do assume that $P$ has been
  replaced by $U^*PU$. But this reduces the norm calculation to the
  direct sum of a set of $2 \times 2$ matrices of the form of the
  lemma together with the diagonal projection $D_4$. Now if $\lambda
  \in \sigma^\prime(C)$, then this $2 \times 2$ matrix is necessarily
  rank one and so the lemma applies. Note also that in this case the
  corresponding eigenvalue of $D_1$ is $1-\lambda$ and that
  $|1-2(1-\lambda)|=|-1+2\lambda|=|1-2\lambda|$ so the values of this
  function agree. When $\lambda=0,$ then this $2
  \times 2$ matrix is either the 0 matrix or it is rank 1 and the
  corresponding eigenvalue of $D_1$ is $1$.
\end{proof}

We now provide a counterexample to Conjecture A.

\begin{thm} Let $\{f_1,...,f_n \}$ be a uniform equiangular Parseval frame for
  $\bb C^k$ with
  $n>2k$ and let $P=(\langle f_i,f_j \rangle )$ be the correlation
  matrix. If there exists a diagonal symmetry, S, such that,
  $\|PSP\|\le 2\delta_P= \frac{2k}{n},$ then $(k-1)n^2 \le 4k^2(n-1).$
\end{thm}
\begin{proof} Without loss of generality we may assume that $S$ is a diagonal symmetry with $m$ diagonal entries
  that are $+1$ and $n-m$ diagonal entries that are $-1$ and, $m \le
  n-m.$
Putting $P$ into the form of the Proposition, we see that since $D_4$
  is a projection, if it is non-zero, then $\|PSP\|=1$. So we
  may assume that $D_4=0.$

Similarly, if any of the diagonal entries of $D_1$ or $D_3$ are 1,
 then $\|PSP\|=1.$ Thus, when we put $P$ into the form of the above Proposition, we obtain a direct sum of $2 \times 2$ rank 1 projections, together with some matrices of all 0's.

Let $0 < \lambda_1 \le ... \le \lambda_t < 1,$ denote the non-zero diagonal entries of $D_1$, so that the corresponding diagonal entries of $D_3$ are $1-\lambda_1,..., 1-\lambda_t,$ and the remaining entries of $D_3$ are 0's.  By the above Theorem, we have that $\|PSP\| = \max \{ |1-2\lambda_1|, |1-2\lambda_k| \}= \max \{ 1- 2\lambda_1, 2\lambda_k -1 \} \le \frac{2k}{n}.$
Hence, $\frac{n-2k}{2n} \le \lambda_1$ and $\lambda_k \le \frac{n+2k}{2n}.$

Since $P$ is a rank $k$ projection, we have that $k=
 Tr(P)=Tr(D_1)+Tr(D_3) =t.$ 
Since $Tr(D_1)= Tr(A)= mk/n$, we have that $0 <
  \lambda_1 \le m/n \le \lambda_k.$ Hence, $\frac{n-2k}{2n} \le m/n \le \frac{n+2k}{2n}$ yielding $n \le 2k+2m$. Note also, that by the choice of $m$ we have that $2m \le n,$ so that the other inequality is automatically satisfied.

If we let, $\mu_1,..., \mu_k$ be the corresponding entries of $D_2$, then since each matrix, $\begin{pmatrix} \lambda_i & \mu_i \\ \mu_i & 1-\lambda_i \end{pmatrix}$ is a rank one projection and since $\mu_i \ge 0,$ we have that $\mu_i^2 = \lambda_i(1-\lambda_i).$

Since $P$ is the correlation matrix of a uniform equiangular (n,k)-frame, by \cite{HP}, we have that every off-diagonal entry of $P$ is of constant modulus, $c= \sqrt{\frac{k(n-k)}{n^2(n-1)}}.$ This yields,
$$\sum_{i=1}^k \mu_i^2 = Tr(B^*B)= m(n-m)c^2 \le \frac{n^2c^2}{4} = \frac{k(n-k)}{4(n-1)}.$$

Now observe that the function $t(1-t)$ is increasing on [0,1/2] and decreasing on [1/2,1]. Thus, we have that $\min \{ \lambda_1(1- \lambda_1), \lambda_k(1-\lambda_k) \} = \min \{ \mu^2_1,..., \mu^2_k \} \le Tr(B^*B)/k \le \frac{n-k}{4(n-1)}.$

However, since $\frac{n-2k}{2n} \le \lambda_1,$ we have $\frac{n-2k}{2n}(1 - \frac{n-2k}{2n})= \frac{n^2 -4k^2}{4n^2} \le \lambda_1(1-\lambda_1).$  Similarly, using the fact that $1/2 < \frac{n+2k}{2n},$ one sees that $\frac{n+2k}{2n}( 1- \frac{n+2k}{2n}) = \frac{n^2-4k^2}{4n^2} \le \lambda_k(1 - \lambda_k).$

Combining these inequalities, yields $\frac{n^2 -4k^2}{4n^2} \le \frac{n-k}{4(n-1)}.$ Cross-multiplying and canceling like terms yields the result.
\end{proof}

Note that the above inequality, for n and k large becomes asymptotically, $n \le 4k.$ Thus, any uniform, equiangular (n,k)-frame with $n/k >> 4,$ and $n$ sufficiently large will yield a counterexample.

\begin{cor} There exist uniform, equiangular Parseval frames whose projection matrices are counterexamples to Conjecture A.
\end{cor}
\begin{proof} In \cite[Example~6.4]{BP} a real uniform, equiangular (276, 23)-frame
  is exhibited and these values satisfy $(k-1)n^2 > 4 k^2(n-1).$ In \cite{K}, uniform, equiangular (n,k)-frames are constructed using Singer difference sets of size, $$n= \frac{q^{m+1} - 1}{q - 1}, k= \frac{q^m - 1}{q - 1},$$
where $q= p^r$ with $p$ a prime. Note that $n/k > q-1.$ Since Singer
difference sets are known to exist for infinitely large $q$, these frames give a whole family of counterexamples.
\end{proof}

We now turn our attention to Conjecture B. We let $\gamma, \epsilon > 0$ be as in the statement of the conjecture. For each partition of $\{ 1,...,n \}= R \cup T$ into two disjoint sets, $R,T$, we let $Q_R, Q_T$ denote the corresponding diagonal projections.

\begin{thm} \label{Pa}
Let $\gamma, \epsilon >0$ be fixed, let $\{f_1,...,f_n \}$ be a uniform Parseval frame for
  $\bb R^k$ with
  $k/n < \min \{ \gamma, \epsilon/2, 1/2 \}$ and let $P=(\langle f_i,f_j \rangle )$ be the correlation
  matrix.  If Conjecture B is true for the pair $(\gamma, \epsilon)$, then there exists a partition $\{1,...,n\}= R \cup T$ such that $Tr(Q_RPQ_TPQ_R) \ge \frac{k\epsilon(2 - \epsilon)}{4}.$
\end{thm}
\begin{proof} Each such partition defines a diagonal symmetry as before and corresponding to such a partition we write $\begin{pmatrix} A &B\\B^* &C \end{pmatrix}$. Note that $Q_RPQ_T = \begin{pmatrix} 0 & B\\0 & 0 \end{pmatrix}$ so that $Tr(Q_RPQ_TPQ_R)= Tr(BB^*).$

We have that $\delta_P= k/n < \gamma.$ We repeat the proof above, with $m= \min \{ |R|,|T| \}$.

Letting $\lambda_1$ be the minimum non-zero eigenvalue and $\lambda_k$ the largest eigenvalue of $A$ as
before, we have $1- \epsilon \ge \|PSP\| \ge
\max \{ |1-2 \lambda_1|, |1- 2 \lambda_k| \}$ and, hence, $\lambda_1 \ge \epsilon/2$ and $1 - \lambda_k \ge \epsilon/2.$

Using the properties of the function $t \to t(1-t)$ and the fact that $\sum_{i=1}^k \lambda_i(1- \lambda_i) = Tr(B^*B),$ we have that $\epsilon/2(1-\epsilon/2) \le \min \{ \lambda_1(1 - \lambda_1), \lambda_k(1 - \lambda_k) \} \le 1/k Tr(B^*B)$, which yields the result.
\end{proof}

Using equiangular frames we can obtain a relation between $\gamma$ and $\epsilon$ in Conjecture B.

\begin{thm} Assume that Conjecture B is true for a pair $(\gamma, \epsilon)$ and let $\{ f_1,..., f_n \}$ be a uniform, equiangular (n,k)-frame with $k/n \le \gamma.$ Then $\epsilon(2- \epsilon) \le \frac{n-k}{n-1}.$ 
\end{thm}
\begin{proof} By the above theorem, we have that there exists a partition with $|R|=m,$ such that $\frac{k\epsilon(2- \epsilon)}{4} \le  Tr(Q_RPQ_TPQ_R) =  m(n-m) c^2 \le \frac{n^2}{4} c^2= \frac{k(n-k)}{4(n-1)}.$
\end{proof}

If we have that infinitely many uniform, equiangular (n,k)-frames
exist for which $n \to +\infty$ and $k/n \to \gamma,$ then
$$\frac{n-k}{n-1}= \frac{1- k/n}{(1 - 1/n)} \to 1- \gamma,$$ and
hence, $\epsilon(2 - \epsilon) \le 1- \gamma.$  If for a given prime p, there are infinitely many Singer difference sets, with $q=p^r,$ and we choose, $1/q \le \gamma$
then we get that $\epsilon(2 - \epsilon) < \frac{q -1}{q}.$

Unfortunately, there are no
uniform Parseval (n,k) frames which violate the
trace inequality in Theorem \ref{Pa}, so that finding
a counter-example to Conjecture B is more subtle.  We will show this 
below.

First, let us change the notation.  If $\{f_i\}_{i=1}^{n}$ is
a Parseval frame for $\l_2^k$ with analysis operator
$V$ then the frame operator is $S=V^{*}V =I$ and
$P=VV^{*}$ is a projection on $\l_2^n$ onto the image
of the analysis operator (which is now an isometry).
Let $\{R,T\}$ be a partition of $\{1,2,\ldots ,n\}$.
If $x = \sum_{i=1}^{n}a_ie_i$ then
$$
Q_Rx = \sum_{i\in R}a_ie_i.
$$
Next,
$$
PQ_Rx = \sum_{j=1}^{n}\langle \sum_{i\in R}a_if_i,f_j\rangle e_j.
$$
Finally,
$$
Q_T PQ_Rx = \sum_{j\in T}\langle \sum_{i\in R}a_if_i,f_j \rangle e_j.
$$
It follows that
$$
Q_TPQ_Re_i = \sum_{i\in R}\sum_{j\in T}\langle f_i,f_j
\rangle e_j.
$$
Now we have:

\begin{lemma}
Given the conditions above we have
$$
Trace(Q_RPQ_TPQ_R) = \sum_{i\in R}\sum_{j\in T}
|\langle f_i,f_j\rangle |^2.
$$
\end{lemma}

\begin{proof}
We compute:
\begin{eqnarray*}
\sum_{i=1}^{n} \langle Q_RPQ_TPQ_Re_i,e_i \rangle &=&
\sum_{i=1}^{n} \langle Q_TPQ_Re_i,PQ_Re_i \rangle \\
&=&  \sum_{i=1}^{n}\langle Q_TPQ_Re_i,Q_TPQ_Re_i \rangle \\
&=& \sum_{i=1}^{n}\|Q_TPQ_Re_i\|^2\\
&=& \sum_{i\in R}\sum_{j\in T}|\langle f_i,f_j \rangle |^2.
\end{eqnarray*}
\end{proof}

Now we need to recall a result of Berman,
Halpern, Kaftal and Weiss \cite{BHKW}.

\begin{thm}\label{WKHB}
Let $(a_{ij})_{i,j=1}^{n}$ be a self-adjoint matrix with non-negative
entries and with zero diagonal so that
$$
\sum_{m=1}^{n} a_{im} \le B,\ \ \mbox{for all $i = 1,2,\ldots ,n$}. 
$$
Then for every
 $r\in \bb N$ there is a partition $\{A_{j}\}_{j=1}^{r}$ of
$\{1,2,\ldots ,n\}$ so that for every $j=1,2,\ldots ,r$,
\begin{equation}\label{ETE}
\sum_{m\in A_j} a_{im} \le  \sum_{m\in A_{\ell}}a_{im},\ \ 
\mbox{for every $i\in A_j$ and $\ell \not= j$}.
\end{equation} 
\end{thm}

Now we are ready for our result.

\begin{prop}
If $\{f_i\}_{i=1}^{n}$ is a uniform (n,k)-Parseval frame,
then there is a partition $\{R,T\}$ of $\{1,2,\ldots ,n\}$
so that
$$
Trace (Q_RPQ_TPQ_R) \ge \frac{k}{4}(1- \frac{k}{n}).
$$
In particular, if $\frac{k}{n}$ is small then the trace inequality
of Theorem \ref{Pa} holds.
\end{prop}

\begin{proof}
Applying \ref{WKHB} to the matrix of values
$(a_{ij})_{i,j=1}^{n}$ where $a_{ii}=0$ and
$a_{ij}= |\langle f_i,f_j\rangle |^2$ for
$i\not= j$ we can find a partition $\{R,T\}$
of $\{1,2,\ldots ,n\}$ (and without loss of generality we
may assume that $|R|\ge \frac{n}{2}$) satisfying for all
$i\in R$:
$$
\sum_{i\not= j\in R}|\langle f_i,f_j \rangle |^2 \le
\sum_{j\in T}|\langle f_i,f_j \rangle |^2.
$$
It follows that for all $i\in R$:
\begin{eqnarray*}
\frac{k}{n} &=& \sum_{j=1}^{n}|\langle f_i,f_j \rangle |^2\\
&=& \frac{k^2}{n^2} + \sum_{i\not= j\in R}|\langle f_i,f_j \rangle |^2
+ \sum_{j\in T}|\langle f_i,f_j\rangle |^2\\
&\le& \frac{k^2}{n^2} +2 \sum_{j\in T}|\langle f_i,f_j \rangle |^2.
\end{eqnarray*}
It follows that for all $i\in R$
$$
\sum_{j\in T}|\langle f_i,f_j \rangle |^2 \ge \frac{1}{2}
\left ( \frac{k}{n}-\frac{k^2}{n^2}\right ).
$$
Now,
\begin{eqnarray*}
\sum_{i\in R}\sum_{j\in T}|\langle f_i,f_j \rangle |^2 &\ge&
|R|\frac{1}{2}
\left ( \frac{k}{n}-\frac{k^2}{n^2}\right )\\
&\ge& \frac{n}{2}\frac{1}{2}
\left ( \frac{k}{n}-\frac{k^2}{n^2}\right )\\
&=& \frac{k}{4} \left ( 1- \frac{k}{n}\right ).
\end{eqnarray*}
Now, given $0<\epsilon <1$, 
$$
\frac{\epsilon}{2}(1-\frac{\epsilon}{2})
< \frac{1}{4}
$$
So the trace inequaltiy of Theorem \ref{Pa} will hold provided
$$
 \frac{k}{4}(1- \frac{k}{n}) \ge
 \frac{k\epsilon(2 -\epsilon)}{4}),
$$
which is true for $k/n$ small enough.
\end{proof}

In fact, as with the case of equiangular frames, we see that if, $1 - \gamma \ge \epsilon ( 2 - \epsilon),$ then whenever, $\frac{k}{n} \le \gamma,$ we have that
$$Trace(Q_RPQ_TPQ_R) \ge \frac{k}{4}(1- \frac{k}{n}) \ge \frac{k}{4}(1 - \gamma) \ge \frac{k}{4} \epsilon(2 - \epsilon).$$

\section{A Family of Potential Counterexamples}

It is still unknown if the paving conjectures are true even for a smaller 
family of operators known as the Laurent operators.
In this section we introduce a family of Laurent operators that we believe 
are potential counterexamples to the paving conjecture. We also prove some 
results about these operators that lends credence to the belief that 
they might yield counterexamples. For the purposes of this section, 
it will be convenient to replace the countable index set $\bb N$ by $\bb Z.$

Recall that a matrix, $A=(a_{i,j})_{i,j \in \bb Z}$ is called a 
{\em Laurent matrix} if it is constant on diagonals, i.e., 
$a_{i,j}= a(i-j)$ and that in this case $A$ determines a bounded 
operator on $\ell^2(\bb Z)$ if and only if there exists 
$f \in L^{\infty}[0,1]$ such that 
$a(n) = \hat{f}(n) \equiv \int_0^1 f(t) e^{-2\pi int} dt$ 
and in this case we set $A= L_f$ and call it the 
{\em Laurent operator with symbol f.} Indeed, the Laurent operator $L_f$ 
is just the matrix representation of the operator of multiplication by 
$f, M_f$ on the space $L^2[0,1]$ with respect to the orthonormal basis, 
$\{ e^{2\pi int} \}_{n \in \bb Z}.$  So, in particular, $L_f$ 
is self-adjoint with diagonal 0 if and only if $f$ is real-valued a.e. 
and $\int_0^1 f(t) dt =0.$

The problem of paving Laurent operators was first studied 
in \cite{HKW3} where it was shown that Laurent operators with 
Riemann integrable symbols can be paved.
Further work on the relation between Laurent operators and the
Feichtinger conjecture can be found in Bownik and Speegle \cite{BS}.

Note that $L_f$ is a projection if and only if $f= \chi_E$ for 
some measurable set $E$ and $L_f$ is a reflection if and only if $f=2 \chi_E -1,$ for some measurable set $E$. This reflection will have 0 diagonal when $m(E)= 1/2,$ where $m$ denotes Lebesgue measure.
Thus, modulo the change from $\bb N$ to $\bb Z$, the family of 
Laurent operators corresponding to our set $\cl R$ is exactly the 
set of operators of the form, $L_f, f= 2 \chi_E -1, m(E)= 1/2$ and 
to $\cl P_{1/2}$ is the set of operators of the form 
$L_f, f= \chi_E, m(E)=1/2.$

Hence, we are interested in the Laurent operators that arise from 
certain subsets $E$ with $m(E)= 1/2.$ It is known that for every 
$t, 0<t<1,$ there exists a measurable set $E=E_t$ with $m(E)=t,$ and 
such that for every $0<a<b<1, m(E \cap (a,b)) >0$ 
and $m(E^c \cap (a,b)) >0,$ where $E^c= [0,1] \setminus E.$ One way to 
construct such a set is as a countable union of fat Cantor sets.

We believe that the projections and reflections coming from such sets 
for $t=1/2,$ are good candidates for counterexamples to the paving 
conjectures and we outline our reasons below.

\begin{prop} Let $E$ be a set as above for any, $0 <t<1.$ If $f_1,f_2$ 
are continuous functions such that $f_1 \le \chi_E \le f_2, a.e.,$ 
then $f_1 \le 0$ and $1 \le f_2.$
\end{prop}
\begin{proof} Since $\chi_E$ is zero on a set of positive measure in 
every interval, $f_1 \le 0.$ Similarly, $\chi_E$ is one on a set of 
positive measure in every interval and hence, $1 \le f_2.$
\end{proof}

The above inequalities show that $\chi_E$ is far from Riemann
integrable.

\begin{prop} Let $g, h \in L^{\infty}[0,1],$ with $0 \le h \le 1.$ If 
for every $f_1,f_2 \in C[0,1],$ we have that $f_1 \le g \le f_2, a.e.,$ 
implies that $f_1 \le 0, 1 \le f_2,$ then there exists a positive linear 
map, $\phi: L^{\infty}[0,1] \to L^{\infty}[0,1]$ such that $\phi(f) =f$ 
for every $f \in C[0,1]$ and $\phi(g) =h.$ 
\end{prop}
\begin{proof}
First define $\phi$ on the linear span of $C[0,1]$ and $g$ by 
$\phi(f+ \alpha g) =f + \alpha h,$ and note that the inequalities 
imply that if $f+ \alpha g \ge 0,$ then $f + \alpha h \ge 0.$
Hence, $\phi$ is a positive map. Now using the fact that 
$L^{\infty}[0,1]$ is an abelian, injective operator system, 
this map has a (completely) positive extension to all of $L^{\infty}[0,1].$
\end{proof}

 \begin{prop} Let $g, h \in L^{\infty}[0,1],$ with $0 \le h \le 1.$ 
If for every $f_1,f_2 \in C[0,1],$ we have that $f_1 \le g \le f_2, a.e.,$ 
implies that $f_1 \le 0, 1 \le f_2,$ then there exists a completely positive 
linear map, $\phi: B(\ell^2(\bb Z)) \to B(\ell^2(\bb Z))$ such that 
$\phi(L_f) =L_f$ for every Laurent operator with continuous symbol, 
$f \in C[0,1]$ and $\phi(L_g) =L_h.$ 
  \end{prop}
  \begin{proof} The identification of $L^{\infty}[0,1]$ with the space 
of Laurent operators is a complete order isomorphism. Hence, there exists 
a completely positive projection of $B(\ell^2(\bb Z))$ onto the space of 
Laurent operators. The remainder of the proof now follows from the last 
Proposition.
\end{proof}

\begin{thm} Let $E \subset [0,1]$ be a measurable set with $m(E)= 1/2$ such that for every $0<a<b<1, m(E \cap (a,b)) >0$ and $m(E^c \cap (a,b)) > 0$ and let $P$ denote the projection that is the Laurent operator with symbol $\chi_E.$
Then there exist completely positive maps, $\phi, \psi: B(\ell^2(\bb Z)) \to B(\ell^2(\bb Z))$ such that $\phi(L_f) = \psi(L_f) = L_f$ for every Laurent operator with continuous symbol $f,$ but $\phi(P)=0, \psi(P) =I.$
\end{thm}
\begin{proof} Apply the above Proposition with $h=0$ and $h=1,$ respectively.
\end{proof}

Thus, for the Laurent reflection with 0 diagonal, $R= 2P-I,$ we have 
that $\phi(R)= -I, \psi(R) = +I$ even though these maps fix all Laurent 
operators with continuous symbols. In this sense, the "value" of the 
diagonal of $R$ is not very stationary under completely positive maps 
which fix all Laurent operators with continuous symbol. In fact, it 
follows from the theory of completely positive maps, that the maps $\phi$ 
and $\psi$ constructed above are actually bimodule maps over the C*-algebra 
of Laurent operators with continuous symbol. That is,
$\phi(L_{f_1}XL_{f_2}) = L_{f_1}\phi(X)L_{f_2},$ and 
$\psi(L_{f_1}XL_{f_2}) = L_{f_1}\psi(X)L_{f_2}$ for any continuous 
functions, $f_1,f_2$ and any $X \in B(\ell^2(\bb Z)).$

One suspects that the fact that the diagonal of $R$ can be altered so 
dramatically, while fixing so many other operators, might be an 
obstruction to $R$ being paved. It is also intriguing that for a 
suitable choice of the set $E$, one can actually compute the 
coefficients of the Laurent matrix for $R$, albeit as power series.

\section*{Acknowledgments}

The authors would like to thank Don Hadwin for several valuable conservations.

\end{document}